\documentclass{amsart}
\usepackage{amssymb,amsmath,amsfonts,amscd,euscript}
\newcommand{\nc}{\newcommand}

\nc{\slt}{\mathfrak{sl}_2}
\nc{\slth}{\widehat{\slt}}
\nc{\C}{\mathbb C }
\nc{\Z}{\mathbb Z }
\nc{\N}{\mathbb N }
\nc{\al}{\alpha }
\nc{\be}{\beta}
\nc{\La}{\Lambda}
\nc{\la}{\lambda}
\nc{\ve}{\varepsilon}
\nc{\ch}{{\mathop {\rm ch}}}
\nc{\Tr}{{\mathop {\rm Tr}\,}}
\nc{\U}{{\mathop {\rm U}}}
\nc{\bra}{\langle}
\nc{\ket}{\rangle}
\nc{\x}{{\bf x}}
\nc{\pa}{\partial}
\nc{\ld}{\ldots}
\nc{\cd}{\cdots}
\nc{\sm}{\sqrt{2m}\,}
\nc{\hk}{\hookrightarrow}
\nc{\A}{\mathfrak A}
\nc{\qb}[2]{\genfrac{(}{)}{0pt}{}{#1}{#2}_q}

\newtheorem{theo}{Theorem}[section]
\newtheorem{lem}{Lemma}[section]
\newtheorem{prop}{Proposition}[section]
\newtheorem{cor}{Corollary}[section]
\newtheorem{rem}{Remark}[section]
\newtheorem{conj}{Conjecture}[section]

\begin{document}
\author{B.Feigin and E.Feigin}
\title
[Principle subspace for bosonic vertex operators]
{Principle subspace for bosonic vertex operator $\phi_{\sqrt {2m}}(z)$ and
Jack polynomials}

\address{Boris Feigin:
{\it Russia, Chernogolovka 142432, Landau Institute for The\-ore\-ti\-cal
Physics} and 
{\it Russia, Moscow, Independent University of Moscow, Bol'shoi Vlas'evskii, 7
}}
\email{feigin@mccme.ru}
\address{Evgeny Feigin:
{\it Russia, Moscow, Independent University of Moscow, Bol'shoi Vlas'evsky,7}
and {\it Russia, Moscow, Moscow State University, Mech-Math Faculty,
Department of
Higher Algebra, Leninskie gori, 1}}
\email{evgfeig@mccme.ru}

\begin{abstract}
Let $\phi_{\sm}(z)=\sum_{n\in\Z} a_nz^{-n-m}$, $m\in\N$ be bosonic vertex
operator, 
$L$ some irreducible representation of the vertex algebra $\A_{(m)}$, associated
with one-dimensional lattice $\Z l$, generated by the vector $l$, $\bra l, l\ket =2m$.
Fix some extremal vector $v\in L$. We study the principal subspace $\C[a_i]_{i\in\Z}\cdot v$
and its finitization $\C[a_i]_{i>N}\cdot v$. We construct their bases and find characters. 
In the case of finitization basis is given in terms of Jack polynomials.
\end{abstract}

\maketitle

\section*{Introduction}
Let $L_{0,1}$ and $L_{1,1}$ be irreducible representations of 
 the Lie algebra $\slth$, $v_{2n}\in L_{0,1}, v_{2n+1}\in L_{1,1}$ the set of extremal 
vectors (for example, $v_0$ is a vacuum vector).  
For $x \in \slt$ consider the current $x(z)= \sum_{i \in \Z} x_i z^{-i-1}$
(here we use the notation $x_i=x \otimes t^i$).
Let $e,h,f$ be standard basis of $\slt$. Then
$e_i v_p=0$ for $i\ge p$. Consider the principle subspace 
$V_p=\C[e_{p-1}, e_{p-2},\ld]\cdot v_p$. Let us list some properties of $V_p$ 
(see \cite{fus,one,two,cp}).

1) $V_p\simeq \C[e_{p-1},e_{p-2},\ld]/I_p$, where $I_p$ is an ideal, generated by
coefficients of series $(e_{p-1}+ze_{p-2}+z^2e_{p-3}+\ld)^2$.
 
2) Elements
$$e_{i_1} \cdots e_{i_k}v,\ i_1<p,\ i_{\al}-i_{\al+1} \ge 2, \
k=0,1, \ldots $$
form the basis of $V_p$.
Using this basis one can write formula for the character of $V_p$ and construct
semi-infinite basis of  $L_{0,1}$ and $L_{1,1}$.

3) Consider finitization: the subspace
$V_p(n) \hookrightarrow V_p,\ V_p(n)= \C [e_{p-1}, \dots,e_{p-n}]
\cdot v_p$. Then $\dim V_p(n)=2^n$. 

 Recall that the current $e(z)$ can be 
realized as
bosonic vertex operator $\phi_{\sqrt {2}}(z)$.
 In this paper we generalize the above results to the case of $\phi_
 {\sqrt {2m}}(z)$ for an arbitrary $m \in \N$
 (see  \cite{jim} for the discussion on this topic).

Let $\A_{(m)}$ be lattice vertex algebra, associated with a one-dimensional lattice
$\Z l$, generated by vector $l$ , $\bra l,l \ket =2m$. Let $L_{(m), i}$ be the set of 
irreducible representations of $\A_{(m)}$ (see \cite{dong}). We have
$$L_{(m),i}=\bigoplus_{n\in\Z} H_{\frac{2nm+i}{\sm}},\ \ 0\le i\le 2m-1.$$ 
Here $H_\la$ is highest weight representation of the Heisenberg algebra 
$H=span\{h_i\}$, $i\in\Z$ with a highest weight vector $|\la\ket$. Fix an action of the
Virasoro algebra on $L_{(m),i}$:
$$L_n=\frac{1}{2}\sum_{i+j=n}: h_i h_j : - \frac{m-1}{\sm} (n+1) h_n.$$
Introduce a notation for the extremal vectors of $L_{(m),i}$:
$v_p=|\frac{-p}{\sm}\ket,\ p\in\Z$.  Note that $a_i v_p=0$ for $i\ge p$. Following
\cite{sto} define the principle subspace 
$$V_{(m),p}=\C[a_{p-1},a_{p-2},\ld]\cdot v_p.$$
First we describe such ideal $I_{(m),p}$ that 
$V_{(m),p}\simeq \C[a_{p-1},a_{p-2},\ld]/I_{(m),p}$. Namely, $I_{(m),p}$ is generated
by the coefficients of series 
$$((a_{p-1}+a_{p-2}z+a_{p-3}z^2+\cd)^{(i)})^2,\ \ 0\le i<m,$$
where superscript $(i)$ stands for the $i$-th derivative.
 Using functional   realization (see \cite{sto}) of the dual
 space $V_{(m),p}^*$ we obtain its character, 
 construct the fermionic realization of $V_{(m),p}$ and find 
its
 monomial basis. Namely we prove that elements
 $$a_{i_1} \cdots a_{i_k}v_p, \
 i_1 \le p-1,\ i_{\al}-i_{\al +1} \ge 2m,\ k=0,1,\ld$$
 form the basis of $V_{(m),p}$.
 We also show that this basis is compatible with the operation of taking
 coinvariants.
 It means that the images of vectors $a_{i_1} \cdots a_{i_k} v_p$ with
 $$i_{\al}-i_{\al +1} \ge 2m,\ i_1 \le p-1,\ i_k>N$$
 form basis of the quotient $V_{(m),p}/span\{a_i V_{(m),p},i \le N \}$.

Now consider the subspace
$$V_{(m),p}(n)\hk V_{(m),p},\
 V_{(m),p}(n)=\C [a_{p-1},a_{p-2},\ldots,a_{p-n}]\cdot v_p.$$
 Define numbers $F_i^{(m)}$ as follows:
 $$F_i^{(m)}=i+1,\ i=0,1,\ldots,m-1;\ F_{i+m}^{(m)}=F_i^{(m)}+
 F_{i+m-1}^{(m)},\ i \ge 0.$$
For example, $F^{(1)}_i=2^i$ and $F^{(2)}_i$ are Fibonacci numbers.
We prove that the dimension of $V_{(m),p}(n)$ equals to $F_n^{(m)}$.
The proof is based on the construction of basis of $V_{(m),p}(n)$ in terms
of Jack polynomials $J_\la(m;\x)$. (Which reduce to Schur polynomials in $\slth$ case). 

Namely, let $h_i, i\in\Z$ be the basis of Heisenberg Lie algebra $H$. Then
each space $H_\la$ can be regarded as the space of
polynomials in infinite number of variables via the identification
($h_i$ with  $i<0$ are generating operators)
$$h_{i_1}\cdots h_{i_k} | \la\ket \mapsto \left(\sqrt{\frac{m}{2}}\right)^k 
p_{-i_1}\cdots p_{-i_k},\ i_{\al}<0,$$
where $p_j$ is the $j$-th power sum.
Combining results from \cite{ty,my} (see also \cite{as})
we obtain  that 
$$a^k_0 v_p=c \cdot J_{((p-(k-1)m-1)^k)}(m;\x),$$
where $c$ is some non-vanishing constant and $((p-(k-1)m-1)^k)$ is the partition
$$\underbrace{(p-(k-1)m-1,\ld,p-(k-1)m-1)}_k.$$

For $p\in\Z$ let $V_{(m),p}^{\frac{-p}{\sm}+k\sm} (n)=V_{(m),p}\cap H_{\frac{-p}{\sm}+k\sm}$.
Note that because of the independence (up to a shift of variables) of $V_{(m),p}$ on $p$
spaces $V_{(m),n}(p)$ and $V_{(m),p}(p)$ are isomorphic. Denote the latter space simply by 
$V_{(m)}(p)$.
We show that
$$V_{(m)}^{\frac{-p}{\sm}+k\sm}(p)=\C[h_1,h_2,\ldots] \cdot
a^k_0 v_p$$
and in addition $\C [h_1,h_2,\ldots]\cdot J_{((p-m(k-1)-1)^k)}(m;\x)$ is a
linear span of Jack polynomials $J_{\la}(m;\x)$ with Young diagram
$\la$ being the subdiagram of $((p-m(k-1)-1)^k)$.
This gives us the dimension and character of $V_{(m)}(p)$.

Recall that in \cite{one} spaces $V_{(1)}(n)$ were obtained via the fusion
 procedure. 
Let $Z=(z_1,\ld,z_n)\in\C^n$, $z_i\ne z_j$. Consider the ring
$R=\C[y_1,\ld,y_n]/\bra y_i^2 \ket_{1\le i \le n}$ and pick its generators
$a_i=\sum_{j=1}^n y_j z_j^{i}$, $i=0,\ld,n-1$.
We have $$R\simeq \C[a_0,\ld,a_{n-1}]/I_Z(n),$$
where $I_Z(n)$ is some ideal. Then the following is true:
$$V_{(1)}(n)\simeq \C[a_{0},\ld,a_{n-1}]/ \lim_{Z\to 0} I_Z(n).$$
Conjecturally, the analogous construction exists for the general $m$. Namely one must put
$R$ to be $\C[y_1,\ld,y_n]/\bra y_iy_j \ket_{|i-j|<m}$. Note that the dimension of this 
ring obviously equals $F^{(m)}_n$. We discuss this conjecture in the last section. 

As it was mentioned above $L_{(m),i}$ generalize irreducible representations $L_{0,1}$ and 
$L_{1,1}$ of $\slth$ on the level $1$. Using the semi-infinite 
construction of $L_{(m),i}$ (see \cite{sto} for the case $m=1$) we obtain bases of the latter
spaces.

The paper is organized as follows.

In the first section we consider the quotient of
$\C[\xi_0,\xi_1,\ld]$  by the ideal, generated by  coefficients
of series $((\xi_0+\xi_1z+\ld)^{(k)})^2$, $k=0,1,\ld, m-1$.
We find basis of this quotient and basis of the space of coinvariants.

In the second section we study the principle subspaces $V_{(m),p}$ and
finitizations $V_{(m)}(p)$. We find relations between Fourier coefficients
of $\phi_{\sm}(z)$ and construct
bases of $V_{(m),p}$ and $V_{(m)}(p)$.

Third section is devoted to the computation of character 
and dimension of $V_{(m)}(p)$.

Fourth section contains construction of semi-infinite basis for $L_{(m),i}$.

And in the last section we discuss some possible generalizations and open questions.

{\bf Acknowledgment.} EF wants to thank V.Dotsenko for helpful information
about Jack polynomials. First named author was partially supported by grants
RFBR 02.01.01015, SS 2044.2003.2 and INTAS 03-51-3350.
Second named author was partially supported by grant RFBR 03-01-00167.

\section{Commutative algebras}
Let $\xi_0,\xi_1,\ldots$ be some commuting variables.
Consider an algebra
\begin{equation}
\label{A}
A_{(m)}=\C [\xi_0,\xi_1,\ldots]/\langle (\xi(z)^{(i)})^2\rangle_{0\le i\le m-1}.
\end{equation}
Here $\xi(z)=\sum_{j=0}^{\infty} \xi_j z^j$, \
$\xi(z)^{(i)}$ its $i$-th derivative, and the ideal of relations
(the right hand side of $(\ref{A})$) is
generated by the coefficients of $(\xi(z)^{(i)})^2,\ i=0,1,\ldots, m-~1$.
Surely, $A_{(m)}$ is bi-graded:
$$A_{(m)}=\bigoplus_{k,s\ge 0}A_{(m)}^{k,s},\
A_{(m)}^{k,s}=span \{\xi_{i_1}\cdots
\xi_{i_k},\ \sum_{\al=1}^k i_{\al}=s \} .$$
Define $\ch(A_{(m)})=\sum_{k,s \ge 0}z^kq^s \dim A^{k,s}_{(m)}.$

\begin{lem}
\label{char}
$\ch(A_{(m)})=\sum_{k \ge 0} z^k\frac{q^{mk(k-1)}}{(1-q)\cdots (1-q^k)}$
\end{lem}
\begin{proof}
The dual space $(A_{(m)}^{k,s})^*$ can be realized as a subspace of
symmetric polynomials $f(z_1,\ldots,z_k)$ of degree $s$.
Namely for $\theta\in (A_{(m)}^{k,s})^*$ define
$$f_\theta(z_1,\ldots,z_k)=\sum_{i_1,\ldots,i_k\ge 0}
z_1^{i_1}\cdots z_k^{i_k} \theta (\xi_{i_1}\cdots \xi_{i_k}).$$
One can show that $f=f_\theta$ for some $\theta\in (A_{(m)}^{k,s})^*$ if and
only if
\begin{equation}
f(z_1,\ldots,z_k)=\prod_{1 \le i<j\le k}(z_i-z_j)^{2m} g(z_1,\ldots,z_k)
\label{f}
\end{equation}
with some symmetric $g$.
(This is because of the condition $(\xi(z)^{(i)})^2=0$
for $i=0,\ldots,m-1$). Formula $(\ref{f})$ gives us the character of $A_{(m)}$.
\end{proof}

Now we want to construct the monomial basis of $A_{(m)}$.
In order to do that we embed $A_{(m)}$ into the algebra, generated by fermions.
Consider algebra $F$, generated by variables $\psi_{s}(i),\
s=1,\ldots,2m,\ i \in \N\cup \{ 0\}$, subject to the relations
$\psi_{s}(i)\psi_{t}(j)=-\psi_t(j)\psi_{s}(i)$.
Define the currents $\psi_s(z)=\sum_{i\ge 0}\psi_s(i)z^i$.
Let $\tilde \xi_i \in F, \ i\ge 0$, be elements defined by
$$\tilde \xi(z)=\sum_{i \ge 0}\tilde \xi_iz^i=\psi_1(z)\cdots \psi_{2m}(z).$$

\begin{lem}
Algebra, generated by $\tilde\xi_i$, $i\ge 0$
is isomorphic to $A_{(m)}$ via the
identification $\tilde \xi_i\mapsto \xi_i$.
\end{lem}
\begin{proof}
First it is easy to check that $(\tilde \xi(z)^{(i)})^2=0$ for $i=0,\ldots,
m-1$.
We need to prove that these are defining relations.
To do that it is enough to find the set of linearly independent monomials
in variables $\tilde \xi_0,\tilde \xi_1,\ldots$, whose character (with respect to
the grading $\deg_z \tilde \xi_i=1,\ \deg_q \tilde \xi_i=i$) is equal to the
character of $A_{(m)}$.

Consider the set of admissible monomials
$$\tilde \xi_{i_1}\cdots \tilde
\xi_{i_k}, \ 0\le i_1,\ i_{\al +1}-i_{\al} \ge 2m,\
\al =1,\ldots, k-1.$$
We prove that they are linearly independent.
In fact, let some linear combination of admissible monomials vanish:
\begin{equation}
\sum \al_{i_1,\ldots,i_k}\tilde \xi_{i_1}\cdots \tilde \xi_{i_k}=0.
\label{van}
\end{equation}
Pick monomial $\tilde \xi_{i_1^0}\cdots \tilde \xi_{i_k^0}$ with a
nonzero $\al_{i_1^0,\ld,i_k^0}$
such that for any
other $k$-tuple $(i_1,\ld,i_k)$ with $\al_{i_1,\ldots,i_k}\ne 0$ the
following is true: there exists such $l <k$ that
$$i_1^0=i_1,\ldots,\ i_l^0=i_l,\ i_{l+1}^0>i_{l+1}.$$
We state that $\tilde \xi_{i_1^0}\cdots \tilde \xi_{i_k^0}$ can not be written
as a sum of other monomials from $(\ref{van})$.

In fact,
\begin{equation}
\label{prod}
\tilde \xi_{i_1^0}\cdots \tilde \xi_{i_k^0}=
\prod_{s=1}^k \sum_{\beta_1^s+\ldots
+\beta_{2m}^s=i_s^0}\psi_1(\beta_1^s)\cdots \psi_{2m}(\beta_{2m}^s).
\end{equation}
Right hand side of $(\ref{prod})$ contains nonzero
(because of the condition $i_{\al+1}-i_{\al}\ge 2m$) term
\begin{equation}
\prod_{s=1}^k \psi_1\left(\left[\frac{i_s^0}{2m}\right]\right)
\psi_2\left(\left[\frac{i_s^0+1}{2m}\right]\right)\cdots
\psi_{2m}\left(\left[\frac{i_s^0+2m-1}{2m}\right]\right)
\label{mon}
\end{equation}
($[x]$ is a maximal integer not exceeding $x$).
By definition of $(i_1^0,\ldots,i_k^0)$ the term $(\ref{mon})$ appears
in the sum $(\ref{van})$ only from the product
$\tilde \xi_{i_1^0}\cdots \tilde \xi_{i_k^0}$.
Thus $(\ref{van})$ is impossible.

To finish the proof note that the character of admissible monomials
coincide with the character of $A_{(m)}$.
\end{proof}

\begin{cor}
\label{basam}
$\xi_{i_1}\cdots \xi_{i_k},\ i_{\al +1}-i_{\al}\ge 2m$,
form basis of $A_{(m)}$.
\end{cor}

\begin{prop}
\label{coinbase}
Images of the monomials $\xi_{i_1}\cdots \xi_{i_k},\ i_{\al +1}-i_{\al}
\ge 2m,\ i_k \le N$ in the quotient $A_{(m)}/span \{ \xi_i A_{(m)},\ i>N \}$
form the basis of the latter quotient.
\end{prop}
\begin{proof}
Note that the term $(\ref{mon})$ with $i_k^0 \le N$ can not appear
in the sum of the form
\begin{equation}
\xi_{N+1}x_1+\xi_{N+2}x_2+\ldots+\xi_{N+p}x_p,\ x_i \in A_{(m)}.
\label{sum}
\end{equation}
In fact, substituting $\tilde\xi_i$ instead of $\xi_i$ we obtain that 
$(\ref{sum})$ is a linear combination of the products of the
fermions
$$\prod_{s=1}^{2m} \psi_s(\beta_1^s)\cdots \psi_s(\beta_k^s)$$
with a property that there exists such map 
$\sigma: \{1,\ld,2m\}\to \{1,\ld,k\}$ that 
$\beta_{\sigma(1)}^1+\ldots+ \beta_{\sigma(2m)}^{2m}>N$.
But the latter is not true for $(\ref{mon})$.
Thus the sum $(\ref{sum})$ does not contain monomial $(\ref{mon})$.
Proposition is proved.
\end{proof}

\section{Vertex operators and finitization}
Consider the lattice vertex operator algebra $\A_{(m)}$, associated with 
one-\-dimen\-si\-onal
lattice $\Z l$, generated by vector $l$, $\bra l, l\ket =2m$. 
Let $\phi_{\sqrt{2m}} (z)=\sum_{n\in\Z} a_n z^{-n-m}$ be corresponding bosonic
vertex operator. Consider the set $L_{(m),i}$ of irreducible representations of 
$\A_{(m)}$. Then $$L_{(m), i}=\bigoplus_{n\in\Z} H_{\frac{2nm+i}{\sm}}, \ 0\le i\le 2m-1.$$
Here $H_\la$ is highest weight representation of the Heisenberg algebra $H$ with a basis
$h_i$, $i\in\Z$. We denote the highest weight vector by $|\la\ket$. Consider the action of Virasoro
algebra on $L_{(m),i}$:
$$L_n=\frac{1}{2} \sum_{i+j=n}: h_i h_j : -\frac{m-1}{\sm} (n+1) h_n,$$
where $:\ :$ is a normal ordering sign. 

Spaces $L_{(m),i}$ are bi-graded by operators $h_0$ ($z$-degree) and
$L_0$  ($q$-degree). For any subspace $V\hk L_{(m),i}$ define
$$\ch (z,q,V)=\Tr (z^{h_0} q^{L_0}|_V).$$
Denote 
$$v_p=\left |\frac{-p}{\sm}\right \ket \in \bigoplus_{i=0}^{2m-1} L_{(m),i}.$$ 
For the eigenvector $v$ of the operators $h_0$ and $L_0$ we denote the corresponding
eigenvalues by $\deg_z v$ and $\deg_q v$. For example 
$$\deg_z v_p=\frac{-p}{\sm},\ \ \deg_q v_p=\frac{p^2}{4m}+\frac{p(m-1)}{2m}.$$
Introduce the principle
subspace $V_{(m),p}=\C[a_{p-1},a_{p-2},\ld]\cdot v_p$ (note that 
$a_iv_p=0$ for $i\ge p$ because the difference of $q$-degrees of $v_p$ and $v_{p-2m}$ equals
$p$).

In what follows we use the connection between vertex
operators and Jack polynomials from \cite{my} (see also \cite{as}).
Recall the definition of Jack polynomials (see \cite{mac}).
Let $\Lambda$ be the algebra of symmetric polynomials in infinite
number of variables, $p_k\in\Lambda$ --  power-sums,
$p_k=\sum_{i\in\N} x_i^k$. Define the scalar product on
$\Lambda=\C[p_1,p_2,\ldots]$, depending on the coupling constant $\al$:
$$
\bra p_1^{i_1} \cdots p_k^{i_k}, p_1^{i_1}\cdots p_k^{i_k}\ket_\al=
\al^{-\sum_{s=1}^k i_s} \prod_{s=1}^k s^{i_s} i_s!
$$
and the products of power sums form the orthogonal basis of $\Lambda$
with respect to $\bra\cdot,\cdot\ket_\al$.

Now let $\la=(\la_1\ge\la_2\ldots\ge \la_s \ge 0)$ be  partition. 
Let $l(\la)$ be the length
of $\la$, i.e. the number of nonzero
$\la_i$. Young diagram, attached to $\la$,  is the following subset of $\Z^2$:
$\{ (i,j):\ 1\le i\le l(\la),\ 1\le j\le \la_i\}$.
Jack polynomials $J_\la (\al;\x)=J_\la (\al;x_1,x_2,\ldots)\in~\La$ depend
on the partition $\la$ and coupling constant $\al$. They are uniquely
determined
by following properties ($m_\la$ is symmetrization of the monomial
$x_1^{\la_1}\cdots x_s^{\la_s}$):
\begin{gather*}
J_\la (\al,\x)=m_\la+\sum_{\mu<\la} v_{\la,\mu}(\al) m_\mu\\
\bra J_\la (\al,\x), J_\mu (\al,\x)\ket_\al =0 \text{ if } \la\ne \mu,
\end{gather*}
where for two partitions $\la, \mu$ we write $\mu\le \la$ if
$\mu_1+\ldots +\mu_i\le \la_1+\ldots +\la_i$ for any $i$.
Note that in the similar way one can define Jack polynomials in finite
number of variables. In this case  $J_\la(\al; x_1,\ldots,x_N)=0$ if
$N<l(\la)$.

Now fix $m\in\N$ and consider Fock space $H_s=\C[h_{-1},h_{-2},\ldots]\cdot |s\ket$.
Identify $H_s$ with $\La$ in the following way:
$$
I: h_{-i_1}\cdots h_{-i_k} |s\ket \mapsto
\left(\sqrt{\frac{m}{2}}\right)^k p_{i_1}\ldots p_{i_k}.
$$
Then $h_i, i>0$ acts on $\La$ by $i\sqrt{\frac{2}{m}} \frac{\pa}{\pa p_i}$.
Denote by $\widehat J_\la (m;\x)$ vector $I^{-1}(J_\la (m,\x))$.
Following theorem was proved in \cite{my}:

\begin{theo}
\label{main}
For any $p>0$
the following
equality holds in $H_{\frac{-p}{\sm}+k\sm}$:
$$
a_0^k v_p= C\cdot \widehat J_{((p-m(k-1)-1)^k)} (m; \x),
$$
where $C$ is some nonzero constant and  $(N^k)$ is a partition
$(\underbrace{N,\ldots,N}_k)$.
\end{theo}

Recall the principle subspace $V_{(m),p}$. First we study 
$$V_{(m),-m+1}=\C[a_{-m},a_{-m-1},\ld]\cdot v_{1-m}.$$ We denote it 
simply by $V_{(m)}$.
\begin{lem}
\label{rel}
Let $\phi^{+}_{\sm}(z)=a_{-m}+za_{-m-1}+z^2a_{-m-2}+\cd$. Then 
$$
V_{(m)}\simeq \C[a_{-m}, a_{-m-1},\ldots]/
\bra (\phi^+_{\sm} (z)^{(i)})^2\ket_{0\le i \le m-1},
$$
where $\phi^+_{\sm} (z)^{(i)}$ is the $i$-th derivative.
It means that the ideal of relations in $V_{(m)}$ is generated by the 
coefficients of series $(\phi^+_{\sm}(z)^{(i)})^2$, namely by
\begin{gather*}
\sum_{\genfrac{}{}{0pt}{}{\al_1+\al_2=s}{\al_1,\al_2\ge 0}}
a_{-m-\al_1}a_{-m-\al_2} \al_1(\al_1-1)\cdots (\al_1-i+1)
\al_2(\al_2-1)\cdots (\al_2-i+1),\\   s\ge 0,\ \ i=0,\ld,m-1.
\end{gather*}
\end{lem}
\begin{proof}
Recall that
$$
\phi_{\sm}(z)\phi_{\sm}(w)=(z-w)^{2m}:\phi_{\sm}(z)\phi_{\sm}(w):.
$$
Thus because of $a_i v_{1-m}=0$ for $i>-m$ we obtain
$(\phi_{\sm}^+ (z)^{(i)})^2 v_{1-m}=0$, $0\le i\le m-1$.
Because of the lemma $(\ref{char})$ it is enough to prove that the character
of $V_{(m)}$ is greater or equal (in each weight component) then
$$z^{\frac{m-1}{\sm}} q^{\frac{-(m-1)^2}{4m}} 
\sum_{k=0}^{\infty} \frac{z^kq^{mk^2}}{(1-q)(1-q^2)\cdots (1-q^k)}.$$ 
(Note that
$\deg_q v_{1-m}=\frac{-(m-1)^2}{4m} v_{1-m}, \  \deg_z v_{1-m}=\frac{m-1}{\sm} v_{1-m}.$)

Consider the decomposition $V_{(m)}=\bigoplus_{k,s\ge 0} \tilde V_{(m)}^{k,s}$,
$$\tilde V_{(m)}^{k,s}= span \{a_{i_1}\ld a_{i_k} v_{1-m},\ \ \sum_{\al=1}^k i_\al=-s\}.$$
Denote $\tilde V_{(m)}^k=\bigoplus_s \tilde V_{(m)}^{k,s}$. Identify $(\tilde V_{(m)}^{k})^*$
with the subspace of symmetric polynomials in $k$ variables in the
following way:
$$(\tilde V_{(m)}^k)^*\ni \theta\mapsto f_\theta(z_1,\ld,z_k)=
\sum_{i_1,\ld,i_k\ge 0} z_1^{i_1}\cd z_k^{i_k}
\theta (a_{-i_1}\cdots a_{-i_k} v_{1-m}).$$
Because of the relations on $\phi_{\sm}^+(z)$ we obtain that $f_\theta$ is
of the following form:
\begin{equation}
\label{ta}
f_\theta(z_1,\ld,z_k)=\prod_{i=1}^k z_i^m \prod_{1\le i<j \le k}
(z_i-z_j)^{2m} g(z_1,\ld,z_k),
\end{equation}
where $g$ is some symmetric polynomial. Note that character of the right hand
side of $(\ref{ta})$ with $g$ running over all symmetric polynomials in $k$
variables coincides with
$\frac{q^{mk^2}}{(1-q)(1-q^2)\cdots (1-q^k)}$. Thus we only need to prove
that any symmetric $g$ appears in the right hand side of $(\ref{ta})$ for some
$\theta$.

First note that $ \tilde V_{(m)}^{k,mk^2}$ is a subspace of one-dimensional
space, spanned by vector $v_{1-m-2k}$. At the same time
one can show that $\tilde V_{(m)}^{k,mk^2}$ is not trivial.
Thus for some nontrivial element
$\theta\in (\tilde V_{(m)}^{k,mk^2})^*$
$$f_\theta(z_1,\ld,z_k)=\prod_{i=1}^k z_i^m \prod_{1\le i<j \le k}
(z_i-z_j)^{2m}.$$

Because of the relations $[h_i,a_j]=\sm a_{i+j}$ operators $h_i$ with
$i>0$ acts on $\tilde V_{(m)}^k$. Let $h_i^*$ be the dual operators, acting on
$(\tilde V_{(m)}^k)^*$. One can see that
$$f_{h_i^*\theta}(z_1,\ld,z_k)=\sm(z_1^i+\ld +z_k^i)f_\theta(z_1,\ld,z_k).$$
To finish the proof it is enough to mention that  power sums generates
the algebra of symmetric polynomials.
\end{proof}

Following lemma can be easily verified.
\begin{lem}
\label{isomor}
The map $\C[a_{-m},a_{-m-1},\ld]\cdot v_{1-m}\to \C[a_{p-1},a_{p-2},\ld]\cdot v_p$,
$$a_{i_1}\cd a_{i_k} v_{1-m} \mapsto a_{i_1+m+p-1}\cd a_{i_k+m+p-1} v_p$$
is well-defined isomorphism. 
\end{lem}

Now we will study the finitization of $V_{(m),p}$. Namely for $p>0$ define the subspace
$$V_{(m),p}\hookleftarrow V_{(m)}(p)=
\C[a_{p-1},a_{p-2},\ld,a_0]\cdot v_p.$$
Introduce a notation
$$V_{(m)}^{k,s}(p)=\{v\in V_{(m)}(p):\ h_0v=kv, L_0 v=sv\},\ \
  V_{(m)}^k(p)=\bigoplus_s V_{(m)}^{k,s}(p).
$$

Following lemma is an immediate consequence from the lemma $(\ref{isomor})$.
\begin{lem}
\label{indep}
Denote $V_{(m),p}(n)=\C[a_{p-1},\ld,a_{p-n}]\cdot v_p$. Then the map
$$V_{(m),p}(n)\to V_{(m)}(n),\ \ a_{i_1}\cd a_{i_k} v_p\mapsto a_{i_1-p+n}\cd a_{i_k-p+n} v_n$$ 
is an isomorphism.
\end{lem}

\begin{lem}
\label{heis}
$V_{(m)}^{\frac{-p}{\sm}+k\sm}(p)=\C[h_1,h_2,\ld]\cdot a_0^k v_p.$
\end{lem}
\begin{proof}
Note that $[h_i,a_j]=\sm a_{i+j}$. We need to prove that for any
$p-1\ge i_1\ge i_2\ge \ld \ge i_k\ge 0$ we have
$$a_{i_1}\cdots a_{i_k} v_p \in \C[h_1,h_2,\ld]\cdot a_0^k v_p.$$
Let us use the induction on $r$: the number of such $\al$ that
$i_\al\ne 0$.

The case $r=0$ is trivial.
Suppose our lemma is proved for some $p$. Then
\begin{multline*}
h_{i_{r+1}} (a_{i_1}\cdots a_{i_r} a_0^{k-r} v_p)=
\sm (k-r) a_{i_1}\cdots a_{i_r} a_{i_{r+1}} a_0^{k-r-1} v_p+\\
\sm a_0^{k-r}
\left(\sum_{j=1}^r a_{i_1}\cdots a_{i_j+i_{r+1}}\cdots a_{i_r}\right)
v_p.
\end{multline*}
Thus by induction assumption
$$a_{i_1}\cdots a_{i_r} a_{i_{r+1}} a_0^{k-r-1}v_p \in
\C[h_1,h_2,\ld]\cdot a_0^k v_p.$$
Lemma is proved.
\end{proof}

As it was mentioned above if we identify $H_\al$ with the space $\La$ of
symmetric polynomials in infinite number of variables, then $h_i$ with $i>0$
will act on $\La$ (up to a nonzero constant) as $\frac{\pa}{\pa p_i}$.
Thus because of the above lemma and theorem
$(\ref{main})$ it is important to study the space
$$\C\left[\frac{\pa}{\pa p_1}, \frac{\pa}{\pa p_2},\ld\right]\cdot
J_\la (m; \x)$$
for the rectangular diagram $\la$.

Consider the polynomial $g\in \C[x_2,x_3,\ld]$. Denote by $(Lg)(x_1,x_2,\ld)$
the polynomial obtained from $g(x_2,x_3,\ld)$ by substituting $x_i$ instead
of $x_{i+1}$, $i=1,2,\ld$.
\begin{prop}
Define a collection of maps $T_n: \La\to \La$, $n=1,2,\ld$:
$$T_n(f(x_1,x_2,\ld))=
L\left(\left.\frac{\pa^nf}{\pa x_1^n} \right|_{x_1=0}\right).$$
Then two algebras of operators on $\La$ are equal:
$$\C\left[\frac{\pa}{\pa p_1},\frac{\pa}{\pa p_2},\ld\right]=
\C[T_1,T_2,\ld].$$
\end{prop}
\begin{proof}
We prove that for any $n>0$
$$\C\left[\frac{\pa}{\pa p_1},\frac{\pa}{\pa p_2},\ld,\frac{\pa}{\pa p_n}
\right]=\C[T_1,T_2,\ld,T_n].$$
Let $f\in\La,f=f(p_1,p_2,\ld)$.
Then there exist polynomials $b_{j,n}(k_1,\ld,k_j),\ j=1,\ld,n,\ k_j\in\N$,
such that
\begin{gather}
\label{2sum}
\frac{\pa^nf}{\pa x_1^n}=\sum^n_{j=1}
\sum_{\genfrac{}{}{0pt}{}{k_1,\ld,k_j\ge 1}
{k_1+\cd+k_j\ge n}}\frac{\pa^jf}
{\pa p_{k_1}\cd\pa p_{k_j}}b_{j,n}(k_1,\ld,k_j)
x_1^{k_1+\cd+k_j-n},\\
\nonumber b_{1,n}(k_1)=k_1(k_1-1)\cd(k_1-n+1).
\end{gather}
We prove $(\ref{2sum})$ by induction on $n$.
For $n=1$ we have
$$\frac{\pa f}{\pa x_1}=\sum_{k\ge 1}\frac{\pa f}{\pa p_k}
\frac{\pa p_k}{\pa x_1}=\sum_{k\ge 1}\frac{\pa f}{\pa p_k}kx_1^{k-1}.$$
Now suppose $(\ref{2sum})$ is true for $n$. Note that
\begin{multline*}
\frac{\pa}{\pa x_1}\left(\sum_{k\ge n}\frac{\pa f}{\pa p_k}k(k-1)\cd(k-n+1)
x_1^{k-n}\right)=\\
=\sum_{k\ge n+1}\frac{\pa f}{\pa p_k}k(k-1)\cd(k-n)x_1^{k-n-1}+\\
+\sum_{\genfrac{}{}{0pt}{}{k_1,k_2\ge 1}
{k_1+k_2\ge n+1}}\frac{\pa^2f}{\pa p_{k_1}
\pa p_{k_2}}k_2k_1(k_1-1)\cd(k_1-n+1)x_1^{k_1+k_2-n-1}.
\end{multline*}
In addition
\begin{multline}
\frac{\pa}{\pa x_1}\left(
\sum_{\genfrac{}{}{0pt}{}{k_1+\cd+k_j\ge n}{k_1,\ld,k_j\ge 1}}
\frac{\pa^jf}{\pa p_{k_1}
\cd\pa_{k_j}}b_{j,n}(k_1,\ld,k_j)x_1^{k_1+\cd+k_j-n}\right)=
\\
\nonumber
\sum_{\genfrac{}{}{0pt}{}{k_1+\cd+k_j\ge n+1}{k_1,\ld,k_j\ge 1}}
\frac{\pa^jf}{\pa p_{k_1}\cd\pa p_{k_j}}
b_{j,n}(k_1,\ld,k_j)(k_1+\cd+k_j-n)x_1^{k_1+\cd+k_j-n-1}+\\
\nonumber
\sum_{\genfrac{}{}{0pt}{}{k_1+\cd+k_j\ge n}{k_1,\ld,k_{j+1}\ge 1}}
\frac{\pa^{j+1}f}{\pa p_{k_1}\cd\pa p_{k_{j+1}}}
k_{j+1}b_{j,n}(k_1,\ld,k_j)x_1^{k_1+\cd+k_{j+1}-n-1}.
\end{multline}
Note that the sum in the last line ranges over 
$k_1,\ld,k_{j+1}$ with $k_1+\ld +k_{j+1}\ge n+1$.
Thus $(\ref{2sum})$ is checked.

Now consider
$L\left(\left. \frac{\pa^nf}{\pa x_1^n} \right|_{x_1=0}\right).$
Because of $L\left(\left. p_i\right|_{x_i=0}\right)=p_i$ we obtain
\begin{multline}
T_n(f)=L\left(\left.\frac{\pa^nf}{\pa x_1^n}\right|_{x_1=0}\right)=\\
=\frac{\pa f}{\pa p_n}n!+\sum^n_{j=2}\ \sum_{k_1+\cd+k_j=n}
\frac{\pa^jf}{\pa p_{k_1}\cd\pa p_{k_j}}b_{j,n}(k_1,\ld,k_j).
\label{tn}
\end{multline}
Now we can prove by induction on $n$ that
\begin{equation}
\C\left[\frac{\pa}{\pa p_1},
\ld,\frac{\pa}{\pa p_n}\right]=\C[T_1,\ld,T_n].
\label{ind}
\end{equation}
For $n=1$ we have
$T_1f=L\left(\left.\frac{\pa f}{\pa x_1}\right|_{x_1=0}\right)=
\frac{\pa f}{\pa p_1}$.
In addition if we know $(\ref{ind})$ for $n$, $(\ref{tn})$  gives us the
statement for $n+1$.
Proposition is proved.
\end{proof}

\begin{rem}
For any $f\in\La$
$$f(x_1,x_2,x_3,\ld)=f(x_2,x_3,\ld)+\sum_{k=1}^{\infty} \frac{1}{k!}x_1^k
(T_kf)(x_2,x_3,\ld).$$
\end{rem}

Recall the definition of skew Jack polynomials $J_{\la/\mu}(\x)$.
(From now on we suppress the coupling constant $m$ in the notation
for Jack and skew Jack polynomials.)
Define $f^{\la}_{\mu\nu}$ as follows:
$J_{\mu}(\x)J_{\nu}(\x)\bra J_\mu, J_\mu\ket^{-1} \bra J_{\nu}, J_{\nu}\ket^{-1}=
\sum_{\la} \bra J_\la, J_\la \ket^{-1} 
f^{\la}_{\mu \nu}J_{\la} (\x)$.
Here $\la,\mu,\nu$ are some partitions.
Then
$$J_{\la/\mu}(\x)=\sum_{\nu}
f^{\la}_{\mu\nu}J_{\nu}(\x),\
\ J_{\la/\mu}(x_1)=\sum_{\nu} 
f^{\la}_{\mu\nu}J_{\nu}(x_1)$$
(surely the summation in the second sum ranges over $\nu$ with $l(\nu)\le 1$).
For two partitions $\mu, \la$ we write $\mu\subset \la$ if the Young diagram
of $\mu$ is a subset of the Young diagram of $\la$.
Note that $J_{\la/\mu}(\x)=0$ unless $\mu\subset \la$.

Following proposition is proved in \cite{mac}.
\begin{prop}
\label{exp}
$\ $

$a)$  $J_{\la}(x_1,x_2,\ld)=\sum_{\mu}J_{\la/\mu}(x_1)J_{\mu}(x_2,x_3,\ld)$.

$b)$
Suppose that $\mu\subset \la$. Then $J_{\la/\mu}(x_1)=0$ unless there exists
such $i$ that $\la_j=\mu_j$ for $j\ne i$. In the latter case
$J_{\la/\mu}(x_1)=\phi_{\la \mu}x_1^{\la_i-\mu_i}$,  $\phi_{\la \mu}$
is some nonzero constant.
\end{prop}

\begin{rem}
Note that $\phi_{\la,\mu}$ depend on the coupling constant $\al$ and do not vanish
for $\al$ being natural number.
\end{rem}

\begin{cor}
\label{TJ}
$T_k(J_{\la}(\x))=k!\sum_{\mu}\phi_{\la\mu}J_{\mu}(\x)$, where
the sum ranges over such partitions $\mu$ that there exists $i$ with
$\la_i-\mu_i=k,\ \la_j=\mu_j$ if $j\ne i$.
\end{cor}
\begin{prop}
\label{base}
Let $\la=(s^r)=(\overbrace{s,s,\ld,s}^r)$.
Then $\C\left[ \frac{\pa}{\pa p_1},\frac{\pa}{\pa p_2},\ld\right] \cdot
J_{\la}(\x)$
is a linear span of $J_{\mu}(\x)$ with $\mu\subset\la$.
\end{prop}
\begin{proof}
Because of the corollary $(\ref{TJ})$ and proposition $(\ref{exp})$ we only 
need
to prove that  $J_{\mu}(\x)\hk\C[T_1,T_2,\ld]\cdot J_{\la}(\x)$ for any 
$\mu\subset \la$. 
Let $\mu=(\mu_1\ge\ld\ge\mu_r)$, $\mu_1\le s$. We use the decreasing
induction on such $j$ that $\mu_j<s$ but $\mu_{j-1}=s$.

In the case $j-1=r$ we have $\mu=\la$. Now let $j\le r$. Let $\tilde\mu$
be partition with
$$\tilde\mu_j=s,\ \ \tilde\mu_i=\mu_i \text{ for } i\ne j.$$
Then by induction assumption $J_{\tilde\mu}(\x)\in\C[T_1,T_2,\ld]\cdot 
J_\la(\x)$. Apply
$T_{s-\mu_j}$ to $J_{\tilde\mu}(\x)$.
$$\frac{1}{(s-\mu_j)!} T_{s-\mu_j} J_{\tilde\mu}(\x)=
\phi_{\tilde\mu \mu} J_\mu(\x)+
\sum_{\nu} \phi_{\tilde\mu \nu} J_\nu(\x),$$
where the sum in the right hand side ranges over such $\nu$ that there
exists $i>j$ with
$$\mu_i-\nu_i=s-\mu_j,\ \ \nu_\al=\tilde\mu_\al \text{ for } \al\ne i.$$
But for any such $\nu$ we have
$\nu_j=s$ and thus $J_\nu(\x)\in\C[T_1,T_2,\ld]\cdot J_{\la} (\x)$.
This gives us
$J_\mu (\x)\in\C[T_1,T_2,\ld]\cdot J_\la (\x)$
(because $\phi_{\tilde\mu \mu}\ne 0$). Proposition is proved.
\end{proof}

Recall the space $V_{(m)}^{\frac{-p}{\sm}+k\sm}(p)$.
\begin{theo}
\label{Jack}
$V_{(m)}^{\frac{-p}{\sm}+k\sm}(p)$ has the basis of the form 
$$\widehat J_{\mu} (\x),\ \ \mu\subset ((p-(k-1)m-1)^k).$$
\end{theo}
\begin{proof}
By lemma $(\ref{heis})$ we have
$$V_{(m)}^{\frac{-p}{\sm}+k\sm}(p)=\C[h_1,h_2,\ld]\cdot a_0^k v_p.$$
Recall that if we identify some $H_\mu$ with the algebra of symmetric polynomials,  
operators $h_i$, $i>0$
will act (up to a nonzero constant) as $\frac{\pa}{\pa p_i}$. Thus our
theorem is a consequence of the theorem $(\ref{main})$ and proposition
$(\ref{base})$.
\end{proof}

\section{Dimensions and characters}
In this section we find the dimensions and characters of $V_{(m)}(p)$
and of the space of coinvariants.

For partition $\mu$ denote $\deg_q \mu=q^{\sum \mu_i}$.
\begin{lem}
\label{ch}
Let $\ch_{s,r}(q)=\sum_{\mu\subset (s^r)} \deg_q \mu$. Then
$\ch_{s,r}(q)=\qb{s+r}{r}.$
\end{lem}
\begin{proof}
First note that $\ch_{s,1}(q)=1+q+\ld +q^s=\qb{s+1}{1}$. In addition
$$
\ch_{s,r}(q)=\sum_{s\ge \mu_1\ge\ld\ge \mu_r\ge 0}
q^{\sum_{i=1}^r \mu_i}=\sum_{j=0}^s
\sum_{\genfrac{}{}{0pt}{}{\mu_r=j}{s\ge\mu_1\ge\ld\ge \mu_r}}
q^{\sum_{i=1}^r \mu_i}=\sum_{j=0}^s q^{rj} \ch_{s-j,r-1}.
$$
To finish the proof it is enough to note that 
$$
\qb{s+r}{r}=
\sum_{j=0}^s q^{jr}\qb{s+r-j-1}{r-1}.
$$
This finishes the proof of the lemma.
\end{proof}

\begin{prop}
\

$a)$ $\ch (z,q,V_{(m)}(p))=
z^{\deg_z v_p} 
\sum\limits_{k\ge 0}z^{k\sm} q^{\deg_q v_{p-2km}} \qb{p-(m-1)(k-1)}{k}$.

$b)$ Let $F^{(m)}_n$ be natural numbers defined by
$$F^{(m)}_i=i+1,\ i=0,1,\ld,m-1;\ F^{(m)}_{i+m}=F^{(m)}_{i+m-1}+F^{(m)}_i.$$
Then $\dim V_{(m),p}(n)=F^{(m)}_n.$
\label{abc}
\end{prop}
\begin{proof}
Consider the element $\widehat J_{\la}(\x)\in H_{\frac{-p}{\sm}+k\sqrt{2m}}$.
Note that
$$L_0 \widehat J_{\la}(\x)=(\deg_q v_{p-2km}+\sum_i\la_i) \widehat J_{\la}(\x).$$
Thus theorem $(\ref{Jack})$ and lemma
$(\ref{ch})$ give us point $a)$.
Note that $b)$ follows from $a)$, because 
\begin{multline*}
\dim V_{(m)}(p+m)=\sum_{k\ge 0} \binom{p+m-(m-1)(k-1)}{k}=\\
\sum_{k\ge 0}\binom{p+m-1-(m-1)(k-1)}{k} +\sum_{k> 0} \binom{p+m-1-(m-1)(k-1)}{k-1}=\\
\dim V_{(m)}(p+m-1)+\sum_{k\ge 0} \binom{p-(m-1)(k-1)}{k}\\ =\dim V_{(m)}(p+m-1)+\dim V_{(m)}(p).
\end{multline*}

Proposition is proved.
\end{proof}

\begin{rem}
There is an obvious embedding $V_{(m),p}(p-1)\hookrightarrow
V_{(m)}(p)$. 
Consider the quotient $V_{(m)}(p)/V_{(m),p}(p-1)$.
Natural way to prove the formula for the dimension of $V_{(m)}(p)$  is to  identify  
above quotient 
with the 
space $V_{(m)}(p-m)$. This was done in \cite{one} for $m=1$, but the straightforward 
generalization failed for $m>1$.
\end{rem}

In the following proposition we find the dimension and character of the
space of coinvariants for $V_{(m),m-1}$. We denote $V_{(m),m-1}$ simply by $V_{(m)}$.

\begin{prop}
Define $C_{(m)}(n)=V_{(m)}/span \{a_i V_{(m)}, i\le -m-n\}$.\\
$a)$ \ $\dim C_{(m)}(n)=F^{(2m)}_n$.\\
$b)$ \ Note that $C_{(m)}(n)$ inherit bigrading from $V_{(m)}$. Thus
its character is well-defined. We have
\begin{equation}
\label{coinch}
\ch(z,q,C_{(m)}(n))=z^{\deg_z v_{1-m}} q^{\deg_q v_{1-m}}  
\sum_{k\ge 0} z^{k\sm}  q^{mk^2} \qb{n-(2m-1)(k-1)}{k}.
\end{equation}
\end{prop}
\begin{proof}
$a)$ follows from the statement (see corollary (\ref{coinbase})) that the
images of vectors
$$a_{-i_1}\cd a_{-i_k} v_{1-m},\ i_{\al+1}-i_\al\ge 2m,\ i_1\ge m,\
i_k\le n+m-1;\ k=0,1,\ld$$
form a basis of $C_{(m)}(n)$.

To prove $b)$ note that it follows from the corollary $(\ref{coinbase})$ that
for $n\ge 2m$
\begin{multline*}
\ch(z,q,C_{(m)}(n))=\ch(z,q,C_{(m)}(n-1))+\\ z^{\deg_z v_{1-m}} q^{\deg_q v_{1-m}}
   z^{\sm}q^{m+n-1} \ch(z,q,C_{(m)}(n-2m)).
\end{multline*}
It is straightforward to check that this relation holds for the right hand
side of $(\ref{coinch})$.
\end{proof}

\section{Semi-infinite bases}
In this section we construct monomial bases for the spaces $L_{(m),i}$.
\begin{lem}
\label{embed}
Recall algebra $A_{(m)}=\C[\xi_0,\xi_1,\ld]/\bra (\xi(z)^{(i)})^2 \ket_{0\le i\le m-1}$. Consider the subspace 
$A_{(m)}^1\hk A_{(m)}$, $A_{(m)}^1=\xi_0 A_{(m)}$. Then the map 
$$\rho: A_{(m)}\to A_{(m)}^1,\ \ \xi_{i_1}\cd \xi_{i_k}\mapsto 
\xi_0 \xi_{i_1+2m}\cd\xi_{i_k+2m}$$
is well-defined isomorphism.  
\end{lem}
\begin{proof}
Recall the fermionic realization $A_{(m)}\hk F$ from section $1$. Note that
$$\xi_0=\psi_1(0)\cd \psi_{2m}(0).$$
Thus we obtain that  $\xi_0 \xi_i=0$ for $i=0,\ld,2m-1$ and $\rho$ is an isomorphism. 
\end{proof}

\begin{lem}
\label{Am}
Consider a map $V_{(m),p}\to A_{(m)}$:
$$a_{j_1}\cd a_{j_s} v_p \mapsto
\xi_{-j_1+p-1}\cd \xi_{-j_s+p-1}.$$
Then this map is well-defined isomorphism.
\end{lem}
\begin{proof}
Follows from lemmas $(\ref{rel})$ and $(\ref{isomor})$.
\end{proof}

\begin{cor}
$L_{(m),i}$ has a basis labeled by the sequences $(s_j)_{j\in\Z}$,
$s_j\in\N\cup \{0\}$ with the following properties:
\begin{enumerate}
\item
There exists such $n\in\Z$  that 
$s_{-i+2mk-1}=1$ for $k\ge n$ and $s_j=0$ if ($j\ge -i+2mn-1$ and $j+i+1$ is not
divisible on $2m$).
\item
There exists such $N\in\Z$ that $s_j=0$ for $j<N$.
\item For any $j\in\Z$ \ $s_j+\ld +s_{j+2m-1}\le 1$.
\end{enumerate}
\end{cor}

\begin{proof}
Rescale vectors $v_p=|\frac{-p}{\sm}\ket$ to satisfy equations
$a_{p-1} v_p=v_{p-2m}$. 
Attach to each sequence $(s_j)_{j\in\Z}$, satisfying properties $(1),(2),(3)$,
the following
vector ($n$ comes from the property $(1)$):
\begin{equation}
\label{at}
\left(\prod_{j< -i+2mn} a_j^{s_j}\right)
v_{-i+2nm} \in L_{(m),i}.
\end{equation}
Then because of $(1)$ for any $n_1> n$
$$
\left(\prod_{j< -i+2n_1m} a_j^{s_j} \right)
v_{-i+2n_1m}=
\left(\prod_{j< -i+2nm} a_j^{s_j} \right)
v_{-i+2nm}.
$$
Note that $L_{(m),i}$ is the limit
$
\lim\limits_{n\to {\infty}} V_{(m), -i+2nm}$:
\begin{equation}
\label{inj}
L_{(m),i}=\ld \hk V_{(m),-i}\hk V_{(m),-i+2m} \hk V_{(m),-i+4m}
\hk\ld
\end{equation}
and each of $V_{(m),i}$ is isomorphic to $A_{(m)}$ in the
sense of lemma $(\ref{Am})$. 
In addition embeddings $(\ref{inj})$ are compatible with the isomorphism from the lemma 
$(\ref{embed})$.
Thus corollary $(\ref{basam})$ gives us
that vectors $(\ref{at})$ form basis of $L_{(m),i}$.
\end{proof}

\section{Discussion}
\subsection{Fusion procedure.}
Recall algebra $A_{(m)}$. Let $A_{(m)}(n)\hk A_{(m)}$ be its finitization:
subalgebra,
generated by $\xi_0,\ld,\xi_{n-1}$. Then
$$A_{(m)}(n)\simeq \C[\xi_0,\ld,\xi_{n-1}]/I_{(m)}(n),$$ where $I_{(m)}(n)$ is
some ideal in $\C[\xi_0,\ld,\xi_{n-1}]$. Now let
$Z=(z_1,\ld,z_n)\in \C^n,\ z_i\ne z_j$ for $i\ne j$.
Let $I_{(m),Z}(n)\hk \C[\xi_0,\ld,\xi_{n-1}]$ be the ideal, generated by the
following elements:
$$(\xi_0+\xi_1z_i+\ld +\xi_{n-1}z_i^{n-1})(\xi_0+\xi_1z_j+\ld
+\xi_{n-1}z_j^{n-1}),\
1\le i\le j\le n,\ |i-j|<m.$$
It was proved in \cite{one} that if
$Z(\ve)=(z_1(\ve),\ld,z_n(\ve))\in\C^n$ is
a family of points of $\C^n$ with the properties $z_i(\ve)\ne z_j(\ve)$
($i\ne j$)
and $\lim_{\ve\to 0} z_i(\ve)=0$ ($1\le i\le n$),
then   $\lim_{\ve\to 0} I_{(1),Z(\ve)} (n) = I_{(1)} (n)$. Following conjecture
generalizes this limit (fusion) procedure to the case of  general $m$.

\begin{conj}
\label{conj}
Let $Z(\ve)=(z_1(\ve),\ld,z_n(\ve))\in\C^n$ be a family with the following
properties:
\begin{enumerate}
\item
$z_i(\ve)\ne z_j(\ve)$ for $i\ne j$ and $\lim_{\ve \to 0} z_i(\ve)=0$,
$1\le i\le n$.
\item
$$\lim_{\ve\to 0} \frac{z_2(\ve)-z_1(\ve)}{z_1(\ve)}=0, \ \
  \lim_{\ve\to 0} \frac{z_{i+1}(\ve)-z_i(\ve)}{z_i(\ve)-z_{i-1}(\ve)}=0, \
  i=2,\ld,n-1.
$$
\end{enumerate}
Then $\lim_{\ve\to 0} I_{(m),Z(\ve)} (n) = I_{(m)} (n)$.
\end{conj}

\subsection{Finite-dimensional approximation.}
Consider the following version of the above fusion procedure. 
Let $R_{(m)}(n)$ be the following ring: 
$$R_{(m)}(n)=\C[y_1,\ld,y_n]/\bra y_i y_j \ket_{|i-j|< m}.$$
For $Z=(z_1,\ld,z_n)\in\C^n$ with pairwise distinct $z_i$ consider the set of generators
of the latter ring:
$$a_i=\sum_{j=1}^n z_j^i y_j,\ \ i=0,\ld,n-1.$$
Then for some  ideal $\tilde I_{(m),Z}(n)\hk \C[a_0,\ld,a_{n-1}]$ 
we have
$$R_{(m)}(n)\simeq \C[a_0,\ld,a_{n-1}]/\tilde  I_{(m),Z}(n).$$
Denote by $\tilde I_{(m)}(n)$ the limit $\lim_{Z\to 0}\tilde I_{(m),Z}(n)$, where 
$Z$ goes to $0$
with the restrictions from the conjecture $\ref{conj}$.
\begin{lem}
The homomorphism of algebras 
\begin{equation}
\label{inv}
\C[\xi_0,\ld,\xi_{n-1}]/I_{(m)}(n)\to \C[a_0,\ld,a_{n-1}]/\tilde I_{(m)}(n), 
\end{equation}
sending $\xi_i$ to $a_{n-1-i}$, is an isomorphism.
\end{lem}  

Denote the right hand side of $(\ref{inv})$ by $W_{(m)}(n)$. 
Above lemma  and conjecture $(\ref{conj})$ allows us to identify 
$W_{(m)}(n)$ with $V_{(m)}(n)$.
Thus we obtain that $L_{(m),i}$ 
can be constructed as limit
of the finite-dimensional spaces: 
$$L_{(m),i}=W_{(m)}(2m-i)\hk W_{(m)}(4m-i)\hk \ld$$
as in the case $m=1$ (see \cite{two}). 

Thus starting from the finite-dimensional spaces, obtained by the fusion procedure, one can 
construct the infinite-dimensional space as an inductive limit and endow this limit with
the structure of the representation of the lattice vertex algebra.

\subsection{Reconstruction of Jack polynomials from the fusion product $V_{(m)}(p)$.}
Let $v\in V^{\frac{-p}{\sm}+k\sm}_{(m)}(p)$ be a basis vector of the form 
$\widehat J_\mu (m;\x)$ for some $\mu$. Note that Jack polynomial $J_\mu(m;\x)$ can
be reconstructed from $v$ using the action of the annihilating Heisenberg operators.

Namely, let $d=\deg_q v-\deg_q v_{p-2mk}$. Pick some $i_1,\ld,i_r>0$ with $i_1+\cd +i_r=d$.
Then
\begin{equation}
\label{pa}
h_{i_1}\cd h_{i_k} v=c(i_1,\ld,i_k) v_{p-2mk}.
\end{equation}
Recall that after the identification of any $H_s$ with the space of symmetric polynomials
operators $h_i$ with $i>0$ act as $\sqrt{\frac{2}{m}} i\frac{\pa}{\pa p_i}$.
Thus $(\ref{pa})$ gives us
$$\frac{\pa}{\pa p_{i_1}}\cd \frac{\pa}{\pa p_{i_k}} J_\mu (m;\x)=
\left(\frac{m}{2}\right)^{k/2}
\frac{c(i_1,\ld,i_k)}{i_1\cd i_k}.$$
So the numbers $c(i_1,\ld,i_k)$ allow us to reconstruct the corresponding Jack polynomial. 
(Note that in the case $m=1$, where everything is proved for 
the fusion products, Jack polynomials reduce to Schur polynomials).

\subsection {The odd case.}
Consider the case of vertex operator $\phi_{\sqrt{2m-1}}(z)$. Let $a_i$ be its Fourier 
components, $a_ia_j=-a_ja_i$ . Let $L$ be some irreducible representation of the corresponding 
vertex algebra,
$v\in L$ some extremal vector. Take such number $p$ that $a_{p-1} v\ne 0$ and  
 $a_k v=0$ for any $k\ge p$. Then one can consider the principle subspace 
$U=\bigwedge (a_{p-1}, a_{p-2},\ld)\cdot v$ and its finitization 
$U(n)=\bigwedge (a_{p-1},\ld, a_{p-n})\cdot v$. The study of $U$ can be provided by 
the same means as in the even case. At the same time one needs some other tools to understand
the structure of $U(n)$. We have a conjecture concerning its dimension.
\begin{conj}
$\dim U(n)=F^{(m)}_n.$
\end{conj}

\end{document}